\documentclass[12pt]{article}

\usepackage{latexsym}
\usepackage{amsfonts}
\usepackage{amssymb}
\usepackage{amsmath}
\usepackage{mathrsfs}
\usepackage{hyperref}

\newcommand{\be}{\begin{equation}}
\newcommand{\ee}{\end{equation}}
\newcommand{\bea}{\begin{eqnarray}}
\newcommand{\eea}{\end{eqnarray}}
\newcommand{\bean}{\begin{eqnarray*}}
\newcommand{\eean}{\end{eqnarray*}}
\newcommand{\brray}{\begin{array}}
\newcommand{\erray}{\end{array}}

\newtheorem{dfn}{Definition}[section]
\newtheorem{thm}[dfn]{Theorem}
\newtheorem{lmma}[dfn]{Lemma}
\newtheorem{ppsn}[dfn]{Proposition}
\newtheorem{crlre}[dfn]{Corollary}
\newtheorem{xmpl}[dfn]{Example}
\newtheorem{rmrk}[dfn]{Remark}

\newcommand{\bdfn}{\begin{dfn}\rm}
\newcommand{\bthm}{\begin{thm}}
\newcommand{\blmma}{\begin{lmma}}
\newcommand{\bppsn}{\begin{ppsn}}
\newcommand{\bcrlre}{\begin{crlre}}
\newcommand{\bxmpl}{\begin{xmpl}}
\newcommand{\brmrk}{\begin{rmrk}\rm}

\newcommand{\edfn}{\end{dfn}}
\newcommand{\ethm}{\end{thm}}
\newcommand{\elmma}{\end{lmma}}
\newcommand{\eppsn}{\end{ppsn}}
\newcommand{\ecrlre}{\end{crlre}}
\newcommand{\exmpl}{\end{xmpl}}
\newcommand{\ermrk}{\end{rmrk}}

\newcommand{\clh}{\mathcal{H}}

\author{S.P. Murugan and S. Sundar}
\title{An essential representation for a product system over a finitely generated subsemigroup of $\mathbb{Z}^{d}$}

\begin{document}
\maketitle
\begin{abstract}
Let $S \subset \mathbb{Z}^{d}$ be a finitely generated subsemigroup. 
Let $E$ be a product system  over  $S$. We show that there exists an infinite dimensional  separable Hilbert space $\clh$
and a semigroup $\alpha:=\{\alpha_x\}_{x \in S}$ of unital normal $*$-endomorphisms of $B(\clh)$ such that $E$
is isomorphic to the product system associated to $\alpha$. 
\end{abstract}
\noindent {\bf AMS Classification No. :} {Primary 46L55; Secondary 46L99.}  \\
{\textbf{Keywords :}} $E_0^{P}$-semigroups, Essential representations, Product systems.

\section{Introduction}
Product systems of Hilbert spaces over $\mathbb{R}_{+}$ were originally invented by Arveson (\cite{Arveson})
 in his classification programme of  $E_0$-semigroups. Later Dinh in \cite{Dinh} and Fowler and Raeburn in \cite{Fowler}
considered product systems over more general semigroups.

 Let $\clh$ be a complex  separable Hilbert space and   denote the algebra  of bounded operators on $\clh$ by $B(\clh)$.  Let $P$ be a discrete countable cancellative semigroup with identity $e$. By an $E^{P}$-semigroup on $B(\clh)$, we mean a family $\alpha:=\{\alpha_{t}\}_{t \in P}$ of normal $*$-endomorphisms of $B(\clh)$ such that
\begin{enumerate}
\item[(1)] for $s \in P$, $\alpha_{s}$ is non-zero,
\item[(2)] for $s,t \in P$, $\alpha_{st}=\alpha_{s} \circ \alpha_{t}$, and
\item[(3)] $\alpha_{e}$ is the identity homomorphism.
\end{enumerate}
An $E^{P}$-semigroup $\alpha:=\{\alpha_{t}\}_{ t \in P}$ on $B(\clh)$ is called an $E_{0}^{P}$-semigroup if for every $t \in P$, $\alpha_{t}(1)=1$. 

Let $\alpha:=\{\alpha_{t}\}_{ t\in P}$ be an $E^{P}$-semigroup on $B(\clh)$. For $t \in P$, set
\[
E(t):=\{T \in B(\clh):\alpha_{t}(A)T=TA ~~\textrm{for all~}A \in B(\clh)\}.
\]
We have the following. 
\begin{enumerate}
\item[(1)]Let $t \in P$ be given. For $S,T \in E(t)$, $T^{*}S$ is a scalar, for it commutes with every element of $B(\clh)$, which we denote by $\langle   S|T\rangle$. Then $\langle  ~| ~\rangle$ defines an inner product on $E(t)$ and makes $E(t)$ a Hilbert space.
\item[(2)] The linear span of $\{ST : S \in E(s), T \in E(t)\}$ is dense in $E(st)$ for every $s,t \in P$.
\item[(3)] The disjoint union of Hilbert spaces $\displaystyle \coprod_{t \in P}E(t)$ has an associative multiplication where the multiplication is given by the following rule
          \[
          (s,S)(t,T):=(st,ST).\]
\end{enumerate}
 In other words, $\displaystyle \coprod_{t \in P}E(t)$ forms a discrete product system of Hilbert spaces over $P$.

It is known that every discrete product system of Hilbert spaces over $P$ is isomorphic to a product system associated to an $E^{P}$-semigroup ( See Remark \ref{regular}). 
The relevant question in Arveson's theory is the following: Is every  product system  over $P$ isomorphic to a product system associated with an $E_{0}^{P}$-semigroup on $B(\clh)$ where $\clh$ is an infinite dimensional separable Hilbert space ? 

So far the question is resolved for the semigroup $\mathbb{N}$ in the discrete case and for $\mathbb{R}_{+}$ in the continuous case. 
See \cite{Arv} and \cite{Skeide} for more details.  The goal of this paper is to settle this question in the affirmative for subsemigroups of $\mathbb{Z}^{d}$ which are finitely generated.
This paper is heavily inspired by \cite{Arv} and the  purpose of this paper is to note down the fact that Arveson's technique in \cite{Arv} works for finitely generated subsemigroups of $\mathbb{Z}^{d}$.  A little bit of notation. 
For us $\mathbb{N}$ denotes the set $\{0,1,2,3,\cdots\}$.
 
The authors would like to thank Prof. Partha Sarathi Chakraborty for his geometric insight which helped us in proving Lemma \ref{the bijection}.

\section{Discrete product systems}
In this section, we recall from \cite{Fowler}, the notion of a discrete product system of Hilbert spaces over a semigroup. What follows in this section is well known and is based on \cite{Fowler},
 the monograph \cite{Arveson} and \cite{Arv}. The authors claim no originality of what follows in this section. We have included this section to make the paper easily readable and self contained.

Let $P$ be a  countable cancellative semigroup containing the identity element $e$. 

\begin{dfn}
\label{discrete product system}
By a discrete product system of Hilbert spaces over $P$, we mean a set $E$ together with a surjection $p:E \to P$ such that 
\begin{enumerate}
\item [(1)] for $t \in P$, $E(t):=p^{-1}(t)$ is a non-zero complex separable Hilbert space,
\item[(2)] there exists an associative multiplication $E \times E \ni (x,y) \to xy \in E$ such that $p(xy)=p(x)p(y)$ for every $x,y \in E$,
\item[(3)] for $s,t \in P$, there exists a unitary $u_{s,t}:E(s) \otimes E(t) \to E_{st}$ such that $u_{s,t}(x\otimes y)=xy$ for $x \in E(s)$ and $y \in E(t)$, and 
\item[(4)] the fibre over the identity element $e$, $E(e)=\mathbb{C}$ and for $s \in P$, the multiplication maps $E(e) \times E(s) \to E(s)$ and $E(s) \times E(e) \to E(s)$ are 
simply scalar multiplication.
\end{enumerate}
\end{dfn}
We suppress the surjection $p$ and simply write a product system $E$ as $ \displaystyle E=\coprod_{t \in P}E(t)$ where for $t \in P$, $E(t)$ is the fibre over $t$. From now on, by a product system over $P$, we mean  a discrete product system of Hilbert spaces over $P$.

Let $\displaystyle E:=\coprod_{t \in P}E(t)$ and $\displaystyle F:=\coprod_{t \in P}F(t)$ be product systems  over $P$. We say $E$ is isomorphic to $F$ if for every $t \in P$, there exists a unitary operator $\theta_{t}:E(t) \to F(t)$ such that  $\theta_{st}(xy)=\theta_{s}(x)\theta_{t}(y)$ for $s,t \in P$ and $(x,y) \in E(s) \times E(t)$.


\begin{dfn}
Let $\displaystyle E:=\coprod_{t \in P}E(t)$ be a  product system  over $P$.  Let $\clh$ be a separable Hilbert space. By a representation of $E$ on $\clh$, we mean a map $\phi:E \to B(\clh)$ such that 
\begin{enumerate}
\item[(1)] for $x, y \in E$, $\phi(xy)=\phi(x)\phi(y)$, and
\item[(2)] for $s \in P$ and $x ,y \in E(s)$, $\phi(y)^{*}\phi(x)= \langle  x|y \rangle$.
\end{enumerate}
The representation $\phi$ is called essential if for every $t \in P$, $\overline{\phi(E(t))\clh}=\clh$.
\end{dfn}
 Let $\displaystyle E:=\coprod_{t \in P}E(t)$ be a  product system  over $P$ and let $\phi:E \to B(\clh)$ be a representation. 
 Then there exists an $E^{P}$-semigroup, $\alpha^{\phi}:=\{\alpha_{t}^{\phi}\}_{ t \in P}$ such that the product system associated to $\alpha^{\phi}$ is $\displaystyle \coprod_{t \in P} \phi(E(t))$ (See Prop. 1.11 of \cite{Fowler}).  Also for $t \in P$, $\alpha_{t}^{\phi}$ is given by the equation
 \begin{equation}
 \label{description}
 \alpha_{t}^{\phi}(A)=\sum_{i=1}^{d}\phi(v_i)A\phi(v_i)^{*}
 \end{equation}
 where $\{v_{i}\}_{i=1}^{d}$ is an orthonormal basis for $E(t)$. If $d$ is infinite, the sum in Eq.\ref{description} is a strongly convergent sum.
  It is clear that $E$ is isomorphic to the product system associated to $\alpha^{\phi}$. Moreover $\alpha_{t}^{\phi}(1)$ is the projection onto the closed subspace $\overline{\phi(E(t))\clh}$. Thus the $E^{P}$-semigroup $\alpha^{\phi}$ is an $E_0^{P}$-semigroup if and only if $\phi$ is essential.

\begin{rmrk}
\label{regular}
It is known that a product system over $P$ is isomorphic to a product system associated to an $E^{P}$-semigroup. This is due to the fact that any product system has a representation on a  separable Hilbert space. (See Lemma 1.10 of \cite{Fowler}).  
\end{rmrk}

For $t,s \in P$, we write $t \geq s$ if there exists $a \in P$ such that $t=sa$. Since $P$ is cancellative, it follows that for $t,s \in P$, if  $t \geq s$ then there exists a unique element in $P$, denoted $s^{-1}t $, such that $t=s(s^{-1}t)$. 

Let $\displaystyle E:=\coprod_{t \in P}E(t)$ be a product system over $P$.
Let $t,s \in P$ be such that $t \geq s$. Choose $a \in P$  such that $t=sa$.  For $v \in E(s)$ and $w \in E(t)$, there exists a unique element denoted $v^{*}w \in E(a)$ such that $ \langle  x|v^{*}w \rangle= \langle  vx|w \rangle $ for every $x \in E(a)$.
Note that 
\begin{equation}
\label{inequality}
||v^{*}w|| \leq ||v||||w||
\end{equation}
 for $v \in E(s)$ and $w \in E(t)$.

\begin{lmma}
\label{star operation}
Let $\displaystyle E:=\coprod_{t \in P}E(t)$ be a product system over $P$. Let $t,s,r \in P$ be such that $t \geq s$. Then for $v \in E(s)$, $w_1 \in E(t)$ and $w_{2} \in E(r)$, $v^{*}(w_1w_2)=(v^{*}w_1)w_{2}$.
\end{lmma}
\textit{Proof.} Let $a \in P$ be such that $t=sa$. Consider elements $v \in E(s)$, $w_{1} \in E(t)$ and $w_{2} \in E(r)$. To show $v^{*}(w_1w_2)=(v^{*}w_{1})w_{2}$, it is enough to show that 
\[ \langle  v^{*}(w_1w_2)|u \rangle= \langle (v^{*}w_{1})w_{2}|u\rangle\] for every $u \in E_{ar}$. Since $\{xy: x \in E(a), y \in E(r)\}$ is total in $E_{ar}$, it suffices to show that for $x \in E(a), y \in E(r)$, 
\[ \langle v^{*}(w_1w_2)|xy\rangle= \langle  (v^{*}w_{1})w_{2}|xy\rangle.\] To that end, let $x \in E(a)$ and $y \in E(r)$ be given. Calculate as follows to find that 
\begin{align*}
 \langle v^{*}(w_{1}w_{2})|xy \rangle&=  \langle   w_1w_2|v(xy) \rangle \\
                                    &= \langle w_1w_{2}|(vx)y \rangle \\
                                    &= \langle   w_1|vx\rangle \langle  w_2|y\rangle \\
                                     &= \langle v^{*}w_{1}|x \rangle \langle   w_{2}|y \rangle \\
                                     &= \langle (v^{*}w_{1})w_{2}|xy\rangle.
\end{align*}
This completes the proof. \hfill $\Box$

We also need the following Lemma whose proof is obtained by merely translating the  proof of Lemma 2.4 of \cite{Arv} to our setup. Thus we omit the proof.

\begin{lmma}
\label{innerprod}
Let $ \displaystyle E:=\coprod_{t \in P}E(t)$ be a product system over $P$. Let $t,s \in P$ be such that $t \geq s$. Suppose that $\{v_{i}\}_{i=1}^{d}$ is an orthonormal basis for $E(s)$. Here $d$ is the dimension of $E(s)$.  Then 
$\displaystyle \sum_{i=1}^{d}||v_{i}^{*}\xi||^{2}=||\xi||^2$ for every $\xi \in E(t)$.
\end{lmma}

\section{Construction of an essential representation}
Fix $d \geq 1$. 
Let $S \subset \mathbb{Z}^{d}$ be a non-zero finitely generated subsemigroup. Then $S-S$ is a subgroup of $\mathbb{Z}^{d}$ and hence isomorphic to $\mathbb{Z}^{m}$ for some $m$. Thus with no loss of generality, we can assume that  $S-S=\mathbb{Z}^{d}$. 
The semigroup $S$ is fixed for the rest of this section. 

\begin{ppsn}
\label{1 dimensional}
Let $\displaystyle E:=\coprod_{s \in S}E(s)$ be a product system over $S$.  Suppose that $E(s)$ is $1$-dimensional for every $s \in S$. Then there exists an $E_{0}^{S}$-semigroup $\alpha:=\{\alpha_{s}\}_{s \in S}$ on $B(\ell^{2}(\mathbb{Z}^{d}))$ such that the product system associated to $\alpha$ is isomorphic to $E$.
\end{ppsn}
\textit{Proof.} For $s \in S$, choose a unit vector $e_{s} \in E(s)$. Then for every $r,s \in S$, there exists a unique scalar, denoted $\omega(r,s) \in \mathbb{T}$, such that $e_re_s=\omega(r,s)e_{r+s}$. The associativity of the multiplication of the product system $E$ implies that $\omega$ is a multiplier on $S$ i.e. for $r,s,t \in S$, \[\omega(r,s)\omega(r+s,t)=\omega(r,s+t)\omega(s,t).\] By Theorem 2.2 of \cite{Laca_multiplier}, it follows that $\omega$ extends to a multiplier on $\mathbb{Z}^{d}$. We denote the extension also by $\omega$. 

Let $\{\delta_{z}:z \in \mathbb{Z}^{d}\}$ be the standard orthonormal basis for $\ell^{2}(\mathbb{Z}^{d})$. For $x \in \mathbb{Z}^{d}$, let $U_{x}$ be the unitary on $\ell^{2}(\mathbb{Z}^{d})$ defined by the equation 
\[
U_{x}(\delta_{z})=\omega(x,z)\delta_{x+z}.\] 
Note that $U_{x}U_{y}=\omega(x,y)U_{x+y}$ for all $x,y \in \mathbb{Z}^{d}$. 

For $s \in S$, let $\alpha_{s}$ be the automorphism of $B(\ell^{2}(\mathbb{Z}^{d}))$ defined by the formula \[\alpha_{s}(A)=U_{s}AU_{s}^{*}.\] Then $\alpha:=\{\alpha_{s}\}_{s \in S}$ is an $E_{0}^{S}$-semigroup on $B(\ell^{2}(\mathbb{Z}^{d}))$. Let $\displaystyle F:=\coprod_{s \in S}F(s)$ be the product system associated to the $E_{0}^{S}$-semigroup $\alpha$.  It is clear that for $s \in S$, $F(s)$ is spanned by $U_{s}$. For $s \in S$, let $\theta_{s}:E(s) \to F(s)$ be the unitary such that $\theta_{s}(e_s)=U_s$.  Now it is immediate that $\displaystyle \theta:=\coprod_{s \in S}\theta_{s}:E \to F$ is an isomorphism of product systems. This completes the proof. \hfill $\Box$

\begin{rmrk}
Suppose $S=\mathbb{Z}^{d}$ and $\displaystyle E:=\coprod_{s \in \mathbb{Z}^{d}}E(s)$ is a product system over $S$. Then for every $s \in \mathbb{Z}^{d}$, $E(s) \otimes E(-s) \cong \mathbb{C}$. This implies that $E(s)$ is $1$-dimensional for every $s \in S$. 
Thus by Proposition \ref{1 dimensional}, it follows that $E$ is isomorphic to a product system associated to an $E_{0}^{S}$-semigroup. 

\end{rmrk}

Hereafter we assume that $S \neq \mathbb{Z}^{d}$.  Let us make a few preliminary observations regarding the semigroup $S$. 
Let $\{e_1,e_2,\cdots,e_{r}\}$ be a set of generators for the semigroup $S$ i.e. $S=\{\sum_{i=1}^{r}m_ie_i: m_i  \in \mathbb{N}\}$. Let $a:=\sum_{i=1}^{r}e_{i}$.  
For $x,y \in \mathbb{Z}^{d}$, we write $x \geq y$ if $x-y \in S$. We use the above notations for the rest of this paper. We have the following archimedean principle. 

\begin{lmma}
\label{archimedean}
Let $x \in \mathbb{Z}^{d}$ be given. Then there exists $n \geq 1$ such that $na \geq x$. As a consequence, we have $-ka \notin S$ for every $k \geq 1$. 
\end{lmma}
\textit{Proof.} Since $S-S=\mathbb{Z}^{d}$, there exists integers $m_1,m_2,\cdots,m_r$ such that $x=\sum_{i=1}^{r}m_{i}e_{i}$. Let $n \geq 1$ be such that $n \geq m_i$ for each $i$. 
Then $na-x=\sum_{i=1}^{r}(n-m_i)e_i \in S$. 

Suppose $-ka \in S$ for some $k \geq 1$. Then $-a=-ka+(k-1)a \in S$. This implies that $-na \in S$ for every $n \geq 1$. Let $x \in \mathbb{Z}^{d}$ be given. Then there exists $n \geq 1$ such that $na \geq -x$ or in other words, $na+x \in S$. Hence $x=(na+x)+(-na) \in S$. This forces that $S=\mathbb{Z}^{d}$ which is a contradiction since we have assumed that $S \neq \mathbb{Z}^{d}$. 

\begin{lmma}
\label{trivial intersection}
The intersection $\bigcap_{n=0}^{\infty}(S+na)=\emptyset$. Moreover the sequence $\{S+na\}_{n \geq 0}$ is a decreasing sequence of subsets of $S$. 
\end{lmma}
\textit{Proof.} Suppose $y \in \bigcap_{n=0}^{\infty}(S+na)$. Then $y-na \in S$ for every $n \geq 0$. By Lemma \ref{archimedean}, there exists $n_0 \geq 1$ such that $n_0a-y \in S$. Note that $-a=(y-(n_0+1)a)+(n_0a-y) \in S$ which is a contradiction to Lemma \ref{archimedean}.  It is clear that $\{S+na\}_{n\geq 0}$ is a decreasing sequence of subsets of $S$. This completes the proof. \hfill $\Box$

\textbf{Notations:} For $k \geq 0$, let $L_{k}:=(S+ka)\backslash(S+(k+1)a)$. Then Lemma \ref{trivial intersection} implies that $\{L_{k}:k \geq 0\}$ is a disjoint family of subsets of $S$ whose union is $S$. Observe that for $k \geq 0$, $ka \in L_{k}$. Also note that for $k \geq 0$, $\displaystyle S+ka=\coprod_{m \geq k}L_{m}$. Since $\displaystyle S=\coprod_{k\geq 0}L_{k}$, for $s \in S$, there exists a unique non-negative integer denoted $n(s)$ such that $s \in L_{n(s)}$.  Note that for $s \in S$, $n(s+a)=n(s)+1$. Also observe that for $s \in S$, $s-n(s)a \in L_0$ and if $s \in L_0$ then $s+ka \in L_{k}$ for $k \geq 0$. For $z \in \mathbb{Z}^{d}$, let $L_{z}=(L_0+z)\cap S$. Note that for $k \geq 0$, $L_{k}=L_{ka}$. We use the above notations throughout this paper. 

Let $\displaystyle E:=\coprod_{s \in S}E(s)$ be a product system over $S$ which is fixed for the rest of this section. We assume that there exists $s \in S$ such that $E(s)$ is not $1$-dimensional. Our goal in this section is to construct an essential representation of $E$ on an infinite dimensional separable Hilbert space. Let $e \in E(a)$ be a unit vector which is fixed for the rest of this section.

Let $\mathcal{V}$ denote the vector subspace of sections of $E$ which are square integrable over $L_{z}$ for every $z \in \mathbb{Z}^{d}$.  More precisely, let $f:S \to E$ be a section. Then $f \in \mathcal{V}$ if and only if for every $z \in \mathbb{Z}^{d}$, 
\[
\sum_{s \in L_{z}}||f(s)||^{2} < \infty.
\]
As is customary, an empty sum equals zero.  Let $f \in \mathcal{V}$ and $k \geq 0$ be given. We say that $f$ is $k$-stable if $f(s+a)=f(s)e$ for $s \geq ka$. Note that if $f$ is $k$-stable then $f$ is $k_1$-stable for $k_1 \geq k$. Let $f \in \mathcal{V}$ be given. We say that $f$ is stable if $f$ is $k$-stable for some $k \geq 0$. Denote the set of stable sections in $\mathcal{V}$ by $\mathcal{S}$. Note that $\mathcal{S}$ is a vector subspace of $\mathcal{V}$.  

Let $f \in \mathcal{V}$. We say that $f$ is eventually zero if there exists $k \geq 0$ such that $f(s)=0$ for $s \geq ka$. Denote the set of eventually zero sections in $\mathcal{V}$ by $\mathcal{N}$. Note that $\mathcal{N} \subset \mathcal{S}$ and $\mathcal{N}$ is a vector subspace of $\mathcal{S}$. 

Let $f,g \in \mathcal{S}$ be given. Since $f$ and $g$ are square integrable over $L_{z}$ for every $z \in \mathbb{Z}^{d}$, it follows that for every $k \geq 0$, the  sum $\displaystyle \sum_{s \in L_k}\langle f(s)|g(s)\rangle$ exists. 

\begin{ppsn}
\label{inner product}
Let $f,g \in \mathcal{S}$.  Then the sequence $\displaystyle \Big(\sum_{s \in L_{k}}\langle f(s)|g(s)\rangle\Big)_{k=1}^{\infty}$ converges.
\end{ppsn}
\textit{Proof.} Without loss of generality, we can assume that $f$ and $g$ are $k_0$-stable for some $k_0 \geq 0$.  Let $k \geq k_0$ be given. Note that the map $L_{k_0} \ni s \to s+(k-k_0)a \in L_{k}$ is a bijection. Now calculate as follows to observe that 
\begin{align*}
\sum_{s \in L_{k}}\langle f(s)|g(s) \rangle & = \sum_{s \in L_{k_0}}\langle f(s+(k-k_0)a)|g(s+(k-k_0)a\rangle \\
                                                                 &= \sum_{s \in L_{k_0}}\langle f(s)e^{k-k_0}|g(s)e^{k-k_0}\rangle ~~(\textrm{Since $f$ and $g$ are $k_0$-stable}) \\
                                                                 &=\sum_{s \in L_{k_0}}\langle f(s)|g(s) \rangle.
\end{align*}
This shows that the sequence $\displaystyle \Big(\sum_{s \in L_{k}}\langle f(s)|g(s)\rangle\Big)_{k=1}^{\infty}$ is eventually constant and hence converges. This completes the proof. \hfill $\Box$

For $f, g \in \mathcal{S}$, let \[\langle f|g \rangle:=\lim_{k \to \infty}\Big(\sum_{s \in L_k}\langle f(s)|g(s) \rangle \Big).\] Then $\langle~|~\rangle$ defines a  semi-definite inner product on $\mathcal{S}$. Let $f \in \mathcal{S}$ be given. Note that $\langle f|f \rangle=0$ if and only if $f \in \mathcal{N}$. It is straightforward to see that if $f \in \mathcal{N}$ then $\langle f|f \rangle=0$. Now let $f \in \mathcal{S}$ be such that $\langle f|f \rangle =0$. Assume that $f$ is $k_0$-stable for some $k_0 \geq 0$. Then the proof of Proposition \ref{inner product} implies that $\displaystyle \sum_{s \in L_{k}}||f(s)||^{2}=0$ for every $k \geq k_0$. This implies that $f$ vanishes on $L_{k}$ for $k \geq k_0$. Hence $f$ vanishes on $\displaystyle \coprod_{k \geq k_0}L_{k}=S+k_0a$. Consequently, we have $f \in \mathcal{N}$. 
Thus $\langle~|~\rangle$ descends to a positive definite  inner product on $\mathcal{S}/\mathcal{N}$ which we still denote by $\langle~|~\rangle$. Let $\clh$ be the completion of the pre-Hilbert space $\mathcal{S}/\mathcal{N}$. 

\begin{rmrk}
\label{inner product1}
Let $f, g \in \mathcal{S}$. Assume that $f$ and $g$ are $k_{0}$-stable for some $k_0 \geq 0$. Then the  proof of Proposition \ref{inner product} shows that  \[\displaystyle \langle f|g \rangle= \sum_{s\in L_{k_0}}\langle f(s)|g(s)\rangle.\]
\end{rmrk}

\begin{ppsn}
The Hilbert space $\clh$ is  separable and is non-zero. 
\end{ppsn}
\textit{Proof.} For $k \geq 0$, let $\displaystyle \clh_{k}:=\bigoplus_{s \in L_{k}}E(s)$. Clearly $\clh_{k}$ is separable for each $k \geq 0$. Fix $k \geq 0$. Let $\xi \in \clh_{k}$ be given. Define a section $\widetilde{\xi}:S \to E$ by the following formula:
\begin{equation*}
 \widetilde{\xi}(s):=\begin{cases}
 \xi(s-n(s)a+ka)e^{n(s)-k} & \mbox{ if
} s \geq ka ,\cr
    & \cr
    0 & \mbox{elsewhere}.
         \end{cases}
\end{equation*}
Note that the above definition makes sense since for $s \geq ka$, $n(s) \geq k$. 
We claim the following. 
\begin{enumerate}
\item[(1)] The section $\widetilde{\xi} \in \mathcal{V}$ and is $k$-stable.
\item[(2)] For $s \in L_{k}$, $\widetilde{\xi}(s)=\xi(s)$. 
\end{enumerate}
Let $z \in \mathbb{Z}^{d}$ be given. Note that $\displaystyle \sum_{s \in L_{z}}||\widetilde{\xi}(s)||^{2}=\displaystyle \sum_{s \in A}||\widetilde{\xi}(s)||^{2}$ where $A:=\{s \in L_{z}: s \geq ka\}$. If $A$ is empty, there is nothing to prove. Suppose that $A$ is non-empty. 
We claim that the map $A \ni s \to s-n(s)a+ka \in L_{k}$ is injective. Suppose $s_1,s_2 \in A$ be such that $s_1-n(s_1)a+ka=s_2-n(s_2)a+ka$. To show $s_1=s_2$, it is enough to prove $n(s_1)=n(s_2)$. Suppose not. Without loss of generality, we can assume that $n(s_2)>n(s_1)$.  Note that $s_2-z=(s_1-z)+(n(s_2)-n(s_1))a \in S+a$ which contradicts the fact that $s_2 \in L_{z} =((S+z)\backslash(S+z+a)) \cap S$. Let $B$ be the image of the map $A \ni s \to s-n(s)a+ka \in L_{k}$. Now calculate as follows to observe that 
\begin{align*}
\sum_{s \in L_{z}}||\widetilde{\xi}(s)||^{2}&=\sum_{s \in A}||\widetilde{\xi}(s)||^{2} \\
                                                               &=\sum_{s \in A}||\xi(s-n(s)a+ka)||^{2} \\
                                                               &=\sum_{s \in B}||\xi(s)||^{2} \\
                                                               &\leq \sum_{s \in L_{k}}||\xi(s)||^{2}\\
                                                               &<\infty.
\end{align*}
This shows that $\widetilde{\xi} \in \mathcal{V}$. Let $s \geq ka$ be given. Calculate as follows to observe that 
\begin{align*}
\widetilde{\xi}(s+a)&=\xi(s+a-n(s+a)a+ka)e^{n(s+a)-k} \\
                              &=\xi(s+a-(n(s)+1)a+ka)e^{n(s)+1-k} ~~(\textrm{Since $n(s+a)=n(s)+1$}) \\
                              &=\xi(s-n(s)a+ka)e^{n(s)-k}e \\
                              &=\widetilde{\xi}(s)e.
       \end{align*}
This proves that $\widetilde{\xi}$ is $k$-stable. This proves $(1)$. Note that for $s \in L_{k}$, $n(s)=k$. Now $(2)$ follows from the definition. Remark \ref{inner product1} together with $(1)$ and $(2)$ implies that the map $\clh_{k} \ni \xi \to \widetilde{\xi}+\mathcal{N} \in \clh$ is an isometry which we denote by $V_{k}$. 

Let $f \in \mathcal{S}$ be given. Assume that $f$ is $k$-stable for some $k \geq 0$. Let $\xi \in \clh_{k}$ be defined by $\xi(s)=f(s)$.  Suppose $s \geq ka$. Note that $s=t+(n(s)-k)a$ where $t=(s-n(s)a)+ka$. Observe that $t \in L_{k}$ and in particular $t \geq ka$. Since $f$ is $k$-stable it follows that \begin{align*}
f(s)&=f(t+(n(s)-k)a)\\
     &=f(t)e^{n(s)-k}\\
     &=f(s-n(s)a+ka)e^{n(s)-k}\\
     &=\xi(s-n(s)a+ka)e^{n(s)-k}\\
     &=\widetilde{\xi}(s).
     \end{align*}
Thus we have shown that $\widetilde{\xi}-f$ is eventually zero.  Consequently $\widetilde{\xi}+\mathcal{N}=f+\mathcal{N}$. Hence $\{f+\mathcal{N}: f \in \mathcal{S}\}= \bigcup_{k=0}^{\infty}V_{k}\clh_{k}$.
This implies that $\bigcup_{k=0}^{\infty}V_{k}\clh_{k}$ is dense in $\clh$. As each $\clh_{k}$ is separable, it follows that $\clh$ is separable. Since each $\clh_k$ is non-zero, it is clear that $\clh$ is non-zero. This completes the proof. \hfill $\Box$.

We need the following two important lemmas before defining a representation of $E$ on $\clh$. 

\begin{lmma}
\label{the bijection}
Let $k \geq 0$ and $b \geq ka$ be given. For every $x \in L_{k}$, the intersection $\{x+ma: m\geq 0\} \cap L_{b}$ is singleton. For $x \in L_{k}$, let $\chi(x) \in L_{b}$ be such that \[\{\chi(x)\}=\{x+ma: m\geq 0\} \cap L_{b}.\] Then the map $L_{k} \ni x \to \chi(x) \in L_{b}$ is a bijection. 
\end{lmma} 
\textit{Proof.} Recall that $L_{b}=(S+b)\backslash(S+b+a)$ and $L_{k}=(S+ka)\backslash(S+(k+1)a)$. Let $x \in L_{k}$ be given. By Lemma \ref{archimedean}, there exists $m \geq 0$ such that $ma-(b-x)=x+ma-b \in S$. Let $m(x)$ be the least non-negative integer such that $x+m(x)a \in S+b$.  Suppose $m(x)=0$. Since $b \geq ka$ and $x \notin S+(k+1)a$, it follows that $x=x+m(x)a \notin S+b+a$. In this case, $x+m(x)a \in L_{b}$. Now suppose $m(x) \geq 1$. Then by definition $x+(m(x)-1)a \notin S+b$. Hence $x+m(x)a \notin S+b+a$. In this case too, $x+m(x)a \in L_{b}$. This proves that the intersection $\{x+ma: m \geq 0\} \cap L_{b}$ is non-empty. 

Suppose $x+ma \in (S+b)\backslash(S+b+a)$. By the definition of $m(x)$, it follows that $m \geq m(x)$. Suppose $m>m(x)$. Write $m=n+m(x)$ with $n \geq 1$.  Observe that  $x+ma=(x+m(x)a)+na \in S+b+na \subset S+b+a$. Hence $x+ma \in S+b+a$ which contradicts the fact that $x+ma \in (S+b) \backslash (S+b+a)$. As a consequence, we have $m=m(x)$. This implies that the intersection $\{x+ma: m \geq 0\} \cap L_{b}$ is singleton.

Let $\chi$ be the map described in the statement of the Lemma. We claim that $\chi$ is $1$-$1$. Let $x_1,x_2 \in L_{k}$ be such that $\chi(x_1)=\chi(x_2)$. Then $x_1+m(x_1)a=x_2+m(x_2)a$. It is enough to prove that $m(x_1)=m(x_2)$. Suppose not. Without loss of generality, we can assume that $m(x_1)<m(x_2)$. Then $x_1=x_2+(m(x_2)-m(x_1))a \in S+ka+(m(x_2)-m(x_1))a$. Since $S+ka+(m(x_2)-m(x_1))a$ is a  subset of $ S+(k+1)a$, it follows that $x_1 \in S+(k+1)a$. This contradicts the fact that $x_1 \in (S+ka)\backslash(S+(k+1)a)$. Hence $\chi$ is $1$-$1$.

We claim that $\chi$ is onto. Let $y \in (S+b)\backslash(S+b+a)$ be given. Since $y \in S+b$ and $b \geq ka$, it follows that $y=y-0.a \in S+ka$. Hence the set $\{m \in \mathbb{N}:y-ma \in S+ka\}$ is non-empty. We claim that $\{m \in \mathbb{N}: y-ma \in S+ka\}$ is bounded. Suppose not. Then there exists a sequence $(m_\ell)$ such that $m_\ell \to \infty$ and $y-m_\ell a \in S+ka$. By Lemma \ref{archimedean}, it follows that there exists $m_0 \geq 0$ such that $m_0a-y +ka \in S$. Choose $\ell$ such that $m_{\ell} > m_0$. Then 
\[
(m_0-m_\ell)a=(m_0a-y+ka)+(y-m_\ell a -ka) \in S
\]
 which is a contradiction to Lemma \ref{archimedean} since $m_0-m_\ell<0$.  This proves that the set $\{m \in \mathbb{N}: y-ma \in S+ka \}$ is bounded. Let $m_0$ be the largest non-negative integer such that $y-m_0a \in S+ka$. 
Then $y-(m_0+1)a \notin S+ka$ or in other words $y-m_0a \notin S+(k+1)a$. Hence $y-m_0a \in (S+ka)\backslash(S+(k+1)a)$. Set $x=y-m_0a$. Then $y=x+m_0a \in \{x+ma: m \geq 0\} \cap L_{b}$. Since the intersection $\{x+ma: m\geq 0\} \cap L_{b}$ is singleton, it follows that $\chi(x)=y$. This proves that $\chi$ is onto. This completes the proof. \hfill $\Box$

\begin{lmma}
\label{crucial Lemma}
Let $f,g \in \mathcal{S}$ be given. Assume that $f$ and $g$ are $k$-stable for some $k \geq 0$. Let $b \in S$ be such that $b \geq ka$. Then 
\[
\langle f|g \rangle = \sum_{s \in L_{b}}\langle f(s)|g(s) \rangle.
\]
\end{lmma}
\textit{Proof.} Let $\chi:L_{k} \to L_{b}$ be the bijection described in Lemma \ref{the bijection}.  For $x \in L_{k}$, let $m(x) \geq 0$ be the unique non-negative integer such that $\chi(x)=x+m(x)a$. Now calculate as follows to observe that 
\begin{align*}
\sum_{s \in L_{b}}\langle f(s)|g(s) \rangle &= \sum_{x \in L_{k}}\langle f(\chi(x))|g(\chi(x))\rangle \\
                                                                  & =\sum_{x \in L_{k}}\langle f(x+m(x)a)|g(x+m(x)a)\rangle \\
                                                                  &=\sum_{x \in L_{k}}\langle f(x)e^{m(x)}|g(x)e^{m(x)}\rangle ~~(\textrm{Since $f$ and $g$ are $k$-stable})\\
                                                                  &=\sum_{x \in L_{k}}\langle f(x)|g(x) \rangle \\
                                                                  &=\langle f|g \rangle ~~(\textrm{by Remark \ref{inner product1}}).
\end{align*}
This completes the proof. \hfill $\Box$

Let $b \in S$ and $v \in E(b)$ be given. For  $f \in \mathcal{S}$, let $\phi_{0}(v)f:S \to E$ be the section defined by 
\begin{equation*}
 (\phi_{0}(v)f)(s):=\begin{cases}
 vf(s-b) & \mbox{ if
} s \geq b \cr
    & \cr
    0 & \mbox{elsewhere}.
         \end{cases}
\end{equation*}
Let $f \in \mathcal{S}$ be given. We leave it to the reader to verify that $\phi_{0}(v)f \in \mathcal{V}$. Suppose that $f$ is $k$-stable. Choose $k_0 \geq 0$ such that $k_0a \geq b$. Set $k_1=k_0+k$. Let $s \in S$ be such that $s \geq k_1a$. Then calculate as follows to observe that 
\begin{align*}
(\phi_0(v)f)(s+a)&=vf(s+a-b)\\
                          &=vf(s-b)e ~~(\textrm{Since $s-b \geq k_1a -b =ka +(k_0a-b) \geq ka$})\\
                          &=(\phi_{0}(v)f)(s)e.
\end{align*}
This proves that $\phi_0(v)f$ is $k_1$-stable. 

\begin{ppsn}
\label{representation}
Let $b \in S$ and $u,v \in E(b)$ be given. Then for $f \in \mathcal{S}$, \[
\langle \phi_0(u)f|\phi_{0}(v)f\rangle=\langle u |v \rangle \langle f | f \rangle.
\]
\end{ppsn}
\textit{Proof.} Let $f \in \mathcal{S}$ be given. Assume that $f$ is $k$-stable for some $k \geq 0$. Choose $k_0 \geq 0$ such that $k_0a \geq b$ and set $k_1=k_0+k$. Then $\phi_{0}(u)f$ and $\phi_{0}(v)f$ are $k_1$-stable.  Now calculate as follows to observe that 
\begin{align*}
\langle \phi_{0}(u)f | \phi_{0}(v)f \rangle &= \sum_{s \in L_{k_1}}\langle \phi_{0}(u)f(s)|\phi_{0}(v)f(s) \rangle ~~(\textrm{by Remark \ref{inner product1}})\\
                                                               &=\sum_{s \in L_{k_1}}\langle  u f(s-b)|v f(s-b) \rangle\\
                                                               &=\langle u|v \rangle \sum_{s \in L_{k_1a}}\langle f(s-b)|f(s-b) \rangle \\
                                                               &=\langle u|v \rangle \sum_{s \in L_{k_1a-b}}\langle f(s)|f(s) \rangle\\
                                                               &=\langle u|v \rangle \langle f|f \rangle. ~~(\textrm{Since $k_1a-b \geq ka$ and by Lemma \ref{crucial Lemma}})
\end{align*}
In the above calculation, to obtain the fourth equality, we have used the fact that the map $L_{k_1a} \ni s \to s-b \in L_{k_1a-b}$ is a bijection. This completes the proof. \hfill $\Box$

Let $b \in S$  and $v \in E(b)$ be given.  Prop. \ref{representation} implies that for $f \in \mathcal{S}$, \[\langle  \phi_{0}(v)f|\phi_{0}(v)f \rangle=||v||^{2}\langle  f|f\rangle.\] As a consequence, it follows that there exists a unique bounded linear operator, denoted $\phi(v)$, on $\clh$ such that 
$\phi(v)(f+\mathcal{N})=\phi_{0}(v)f+\mathcal{N}$ for every $f \in \mathcal{S}$. Prop. \ref{representation} implies that for $u,v \in E(b)$, $\phi(v)^{*}\phi(u)= \langle u|v\rangle$. It is clear that $\phi:E \to B(\clh)$ is multiplicative. Thus $\phi$ is a representation of $E$ on $\clh$.
Our goal is to show that $\phi$ is essential.

\begin{rmrk}
The Hilbert space $\clh$ is infinite dimensional. To see this, observe that we have assumed that there exists $b \in S$ such that $E(b)$ is not $1$-dimensional. Let $\{v_{i}\}_{i=1}^{d}$ be an orthonormal basis for $E(b)$ where $d$ is the dimension of $E(b)$. Since $\phi$ is a representation $\{\phi(v_i)\}_{i=1}^{d}$ is a family of isometries with orthogonal range projections. But $d \geq 2$. This implies that $\clh$ is infinite dimensional. 
\end{rmrk}

Let $v \in E(a)$ and $f \in \mathcal{S}$ be given. Define a section $f_{v}:S \to E$ by the following formula \[f_{v}(s)=v^{*}f(s+a).\] We leave it to the reader to verify that $f_{v} \in \mathcal{S}$. We merely indicate that to show $f_{v} \in \mathcal{V}$, one needs to use Eq. \ref{inequality} and to show that $f_{v}$ is stable one needs to use Lemma \ref{star operation}. Note that if $f$ is $k$-stable then $f_{v}$ is $k$-stable. 

\begin{lmma}
Let $v \in E(a)$ be given. For $f \in \mathcal{S}$, $\phi(v)^{*}(f+\mathcal{N})=f_{v}+\mathcal{N}$.
\end{lmma}
\textit{Proof.} Let $f \in \mathcal{S}$ be given. To show that $\phi(v)^{*}(f+\mathcal{N})=f_{v}+\mathcal{N}$, it suffices to show that for $g \in \mathcal{S}$, $\langle \phi(v)^{*}(f+\mathcal{N})|g+\mathcal{N}\rangle=\langle f_{v}+\mathcal{N}|g+\mathcal{N}\rangle$. Let $g \in \mathcal{S}$ be given.  Without loss of generality, we can assume that $f$ and $g$ are $k$-stable for some $k \geq 0$. Then $\phi_0(v)g$ is $k+1$-stable. Now calculate as follows to observe that 
\begin{align*}
\langle \phi(v)^{*}(f+\mathcal{N})|g+\mathcal{N}\rangle & = \langle f+\mathcal{N}|\phi(v)(g+\mathcal{N})\rangle \\
                                                                                       &=\sum_{s \in L_{k+1}}\langle f(s)|\phi_0(v)g(s)\rangle~~(\textrm{by Remark \ref{inner product1}}) \\
                                                                                       &=\sum_{s \in L_{k+1}}\langle f(s)|vg(s-a) \rangle \\
                                                                                       &=\sum_{s \in L_{k+1}}\langle v^{*}f(s)|g(s-a) \rangle \\
                                                                                       &=\sum_{s \in L_{k}}\langle v^{*}f(s+a)|g(s) \rangle \\
                                                                                       &=\langle f_{v}+\mathcal{N}|g+\mathcal{N}\rangle. ~~(\textrm{by Remark \ref{inner product1}})
\end{align*}
In the above calculation, to obtain the fifth equality, we have used the fact that the map $L_{k} \ni s \to s+a \in L_{k+1}$ is a bijection. This completes the proof. \hfill $\Box$

Recall that $\{e_i: i=1,2,\cdots,r\}$ are the chosen generators of $S$ and $\displaystyle a=\sum_{i=1}^{r}e_{i}$. 
\begin{thm}
The representation $\phi$ is essential.
\end{thm}
\textit{Proof.} Let $\alpha:=\{\alpha_{s}\}_{s \in S}$ be the $E^{S}$-semigroup associated to $\phi$. To show that $\alpha_{s}$ is unital for every $s$, it suffices to show that $\alpha_{a}$ is unital. To see this, note that $\alpha_{s}(1) \leq \alpha_{t}(1)$ if $s \geq t$. 
Hence if $\alpha_{a}$ is unital, it follows that $\alpha_{e_i}$ is unital for every $i=1,2,\cdots,r$. But $S$ is generated by $\{e_{i}:i=1,2,\cdots,r\}$. This forces that $\alpha_{s}$ is unital for every $s \in S$ provided $\alpha_{a}$ is unital. 

Let $\{v_{i}\}_{i=1}^{d}$ be an orthonormal basis for $E(a)$ where $d$ denotes the dimension of $E(a)$. We claim that $\displaystyle \sum_{i=1}^{d}\phi(v_i)\phi(v_i)^{*}=1$. Here the sum is interpreted in the strong sense if $d$ is infinite. 
Since $\{\phi(v_{i})\phi(v_{i})^{*}\}_{i=1}^{d}$ forms a mutually orthogonal family of projections, it is enough to show that 
\[
\sum_{i=1}^{d} \langle \phi(v_{i})\phi(v_{i})^{*}(f+\mathcal{N})|f+\mathcal{N} \rangle=||f+\mathcal{N}||^{2}
\]
for every $f \in \mathcal{S}$. 

Let $f \in \mathcal{S}$ be given. Assume that $f$ is $k$-stable. Then $f_{v_i}$ is $k$-stable for every $i$. Now calculate as follows to observe that 
\begin{align*}
\sum_{i=1}^{d}||\phi(v_i)^{*}(f+\mathcal{N})||^{2} &=\sum_{i=1}^{d}||f_{v_i}+\mathcal{N}||^{2} \\
                                                                            &=\sum_{i=1}^{d}\sum_{s \in L_{k}}||f_{v_i}(s)||^{2}\\\
                                                                            &=\sum_{i=1}^{d}\sum_{s \in L_{k}}||v_{i}^{*}f(s+a)||^{2}\\
                                                                            &=\sum_{s \in L_{k}}\sum_{i=1}^{d}||v_{i}^{*}f(s+a)||^{2} \\
                                                                            &=\sum_{s \in L_{k}}||f(s+a)||^{2} ~~(\textrm{ by Lemma \ref{innerprod}})\\
                                                                            &=\sum_{s \in L_{k}}||f(s)||^{2}~~(\textrm{since $f$ is $k$-stable}) \\
                                                                            &=||f+\mathcal{N}||^{2}~~(\textrm{by Remark \ref{inner product1}}).
\end{align*}
In the fourth equality of the above calculation, we have interchanged the order of summation which is permissible since the terms involved are non-negative. This completes the proof. \hfill $\Box$

\bibliography{references}
 \bibliographystyle{plain}
 
 \noindent{\sc S.P. Murugan} (\texttt{spmurugan@cmi.ac.in})\\
         {\footnotesize  Chennai Mathematical Institute, \\
Siruseri, 603103, Tamilnadu.}\\[1ex]
{\sc S. Sundar}
(\texttt{sundarsobers@gmail.com})\\
         {\footnotesize  Chennai Mathematical Institute,  \\
Siruseri,  603103, Tamilnadu.}

\end{document}